# Paraconsistent Existential Graphs Gamma Peirce System

Manuel Sierra-Aristizabal
[1]*(EAFIT University, Colombia)*

***Abstract:***
*In this paper, the paraconsistent propositional logic LG is presented, along with its semantic characterization. It is shown that LG's set of theorems corresponds to the set of valid existential graphs, GET, which turns out to be an extension of Peirce's Gamma system, without becoming Zeman's gamma-4 system. All evidence is presented in a complete, rigorous, and detailed manner. This result is generalized by constructing the paraconsistent system of existential graphs GET4, and its semantic-deductive characterization. Finally, Zeman's Gamma-4, Gamma-4.2, and Gamma-5 existential graph systems are proven to be paraconsistent.*

## I. Introduction

Existential graphs, alpha, beta, and gamut, were created by Charles Sanders Peirce in the late 19th century, see Roberts (1992) and Peirce (1965). Alpha graphs correspond to classical propositional calculus, beta graphs correspond to classical logic of first-order relations. Gamma charts were introduced by Peirce and later extended by Jay Zeman, constructing existential graphs for modal logics S4, S4.2 and S5 in Zeman (1963). On the other hand, Brade and Trymble (2000) have proposed categorical models for alpha existential graphs. Recently, existential graphs were presented for intuitionistic propositional calculus in Oostra (2010) and Oostra (2021), for intuitionistic relationship calculus in Oostra (2011), and for modal logics S4, S4.2, and S5, intuitionist versions, in Oostra (2012). Finally, Sierra (2021) presents the Gamma-LD system of existential graphs, and Sierra (2022) presents the first system of paraconsistent existential graphs.

In this paper, the paraconsistent propositional logic LG is presented, along with its semantic characterization. It is shown that LG's set of theorems corresponds to the set of valid existential graphs, GET, which turns out to be an extension of Peirce's Gamma system, without becoming Zeman's gamma-4 system. All evidence is presented in a complete, rigorous, and detailed manner. This result is generalized by constructing the paraconsistent system of existential graphs GET4, and its semantic-deductive characterization. Finally, Zeman's Gamma-4, Gamma-4.2, and Gamma-5 existential graph systems are proven to be paraconsistent.

## II. Deductive System

In this section, the deductive system of propositional logic GT is presented, its connections with classical propositional calculus, and some of its theorems.

**Definition 1.** The FT set of **GT** formulas is constructed from a set FA of atomic formulas, from the constant $\lambda$ , the unary connective, weak negation $\{-\}$, and the binary connective, conditional $\{\supset\}$ as follows.
$P\in FA$ implies $P\in FT$. $\lambda\in FT$. $X\in FT$ implies $-X\in FT$. $X,Y\in FT$ implies $X\supset\in FT$.
Classical negation, strong affirmation, weak affirmation, biconditional affirmation are defined as:
a) $\sim X = X\supset-\lambda$. b) $+X = \sim-X$. c) $\otimes X = -\sim X$. d) $X\cup Y = \sim X\supset Y$. e) $X\cap Y = \sim(X\supset\sim Y)$. f) $X\equiv Y = (X\supset Y)\cap(Y\supset X)$.

**Definition 2.** The GT system consists of the axioms (where $X,Y,Z\in FT$):
Ax1. $\lambda$
Ax2. $X\supset(Y\supset X)$
Ax3. $[X\supset(Y\ Z\supset)]\supset[(X\supset Y)\supset(X\supset Z)]$
Ax4. $[(X\supset Y)\supset X]\supset X$
Ax5. $-\lambda\supset Z$
Ax6. $(X\supset-\lambda)\supset-X$
Ax7. $[-(X\supset Y)\supset-\lambda]\supset[(-X\supset-\lambda)\supset(-Y\supset-\lambda)]$
Ax+. If $X\in\{Ax1, ..., Ax6\}$ then $-X\supset-\lambda$ is an axiom.
The only rule of inference is the *modus ponens* Mp: from X and $X\supset Z$ we infer Z.

**Definition 3.** For $X, X_1, ..., X_n \in FT$. X is a *theorem of* GT, *denoted* $X \in TT$, if there is a proof of X from the axioms using the rule Mp, i.e., X is the last row of a finite sequence of lines, in which each of the lines is an axiom, or is inferred from two preceding rows, using the inference rule Mp. The number of lines in the sequence is referenced as the *length* of the X proof. Y is a *theorem (or consequence)* of $\{X_1, ..., X_n\}$, or $\{X_1, ..., X_n\}$ implies Y, if there is a proof of Y, from the axioms and assumptions $\{X_1, ..., X_n\}$.

**Proposition 1.** For $X,Y,X_1,..., X_n \in FT$. If $\{X_1, ..., X_n, X\}$ implies Y in GT, then $\{X_1, ..., X_n\}$ implies $X \supset Y$.

Proof. Axioms 2, 3 and 4, with the single inference rule Mp, determine the calculus for the classical implication CIC Rasiowa (1974), in which the deduction theorem, TD, applies.

**Proposition 2.** For $X,Y \in FT$. $+(X \supset Y) \supset (+X \supset +Y) \in TT$.

Proof. Ax7 and definition.

**Proposition 3.** For $X \in FT$. If $X \in TT$ then $+X \in TT$.

Proof. Suppose $X \in TT$, $+X \in TT$ will be tested, by induction over the length of the proof of X.
Base step. The length of the proof of X is 1, i.e., X is an axiom. If X is one of axioms 1 to 7, Ax+ gives +X. If X is one of the axioms generated by Ax+, it is already of the form +Y.
Induction step. As an inductive hypothesis, if the length of the proof of Y is less than L, then +Y is a theorem. Suppose that the proof of X has length L greater than 1. It follows that X is an axiom or X is a consequence of previous steps using the inference rule Mp. In the first case, proceed as in the base step. In the second case, we have, for some formula Z, proofs of $Z \supset X$ and Z, both of which are shorter in length than L. From the inductive hypothesis we infer $+(Z \supset X), +Z \in TT$. By proposition 2 we have $+(Z \supset X) \supset (+Z \supset +X) \in TT$, applying the rule Mp twice we get that $+X \in TT$.
So, according to the principle of mathematical induction, it has been proved $X \in TT$ implies $+X \in TT$.

**Proposition 4.** For $X,Y \in FT$. GT theorems are: a) $(X \supset \sim Y) \supset (Y \supset \sim X)$. b) $\sim(X \supset X) \supset Y$. c) $X \cup \sim X$. d) $X \supset \sim\sim X$.
e) $\sim\sim X \supset X$. f) $(X \supset Y) \supset (\sim Y \supset \sim X)$. $(\sim Y \supset \sim X) \supset (X \supset Y)$.

Proof part a. Suppose $X \supset (Y \supset -\lambda)$, Y, X. By Mp is derived $Y \supset -\lambda$, again by Mp is inferred $-\lambda$. Applying TD 3 times and using the definition of $\sim$, concludes $(X \supset \sim Y) \supset (Y \supset \sim X)$.
Proof part b. Suppose $\sim(X \supset X)$, i.e. $(X \supset X) \supset -\lambda$, but $X \supset X$ is a theorem of CIC, resulting in $-\lambda$, using Ax5 follows Y. Applying TD concludes $\sim(X \supset X) \supset Y$.
Proof part c. By the principle of identity of the CIC we have $\sim X \supset \sim X$, by the definition of $\cup$ we conclude $X \cup \sim X$.
Proof part d. Suppose X, $X \supset -\lambda$. By Mp we follow $-\lambda$, applying TD 2 times and definition of $\sim$ we conclude $X \supset \sim\sim X$.
Proof part e. Suppose $\sim\sim X$, i.e., $\sim X \supset -\lambda$, by Ax5 we have $-\lambda \supset X$, by CIC we deduce $\sim X \supset X$, i.e., $(X \supset \lambda -) \supset X$, using Ax4 implies X. By TD we conclude $\sim\sim X \supset X$.
Proof part f. parts a, d and e.

**Proposition 5**. Sean $X,Y,Z \in FT$. GT theorems are: a) $X \supset (X \cup Y)$. b) $X \supset (Y \cup X)$.
c) $(X \supset Y) \supset [(Z \supset Y) \supset (\{X \cup Z\} \supset Y)]$

Proof part a. Suppose X, $X \supset -\lambda$, i.e., $\sim X$, by Mp we get $-\lambda$, according to Ax5 we derive Y. Applying TD 2 times concludes $X \supset (\sim X \supset Y)$, i.e., $X \supset (X \cup Y)$.
Proof part b. By part a we conclude $X \supset (\sim X \supset Y)$, using proposition 4, it can be said that $X \supset (\sim Y \supset X)$, i.e., $X \supset (Y \cup X)$.
Proof part c. Suppose $X \supset Y$, $Z \supset Y$, $X \cup Z$, i.e., $\sim X \supset Z$, by CPC we infer $\sim X \supset Y$, by proposition 4 we derive $\sim Y \supset X$, by CIC we infer $\sim Y \supset Y$, i.e. $(Y \supset -\lambda) \supset Y$, by Ax4 we get Y. Applying TD 3 times we get $(X \supset Y) \supset [(Z \supset Y) \supset (\{X \cup Z\} \supset Y)]$.

**Proposition 6.** For $X,Y \in FT$. GT theorems are: a) $(X \cap Y) \supset X$. b) $(X \cap Y) \supset Y$.
c) $(X \supset Y) \supset [(X \supset Z) \supset (X \supset \{Y \cap Z\})]$. d) $X \supset [Y \supset (X \cap Y)]$. e) $+(X \cap Y) \equiv (+X \cap +Y)$.

Proof part a. Suppose $X \cap Y$, i.e., $\sim(X \supset \sim Y)$, so $(X \supset \sim Y) \supset -\lambda$, by Ax5 we have $-\lambda \supset X$, by CIC we infer $(X \supset \sim Y) \supset X$, using Ax4 results X. By TD we conclude $(X \cap Y) \supset X$.
Proof part b. Suppose $X \cap Y$, i.e., $\sim(X \supset \sim Y)$, so that $(X \supset \sim Y) \supset -\lambda$, using proposition 4 we deduce $(Y \supset \sim X) \supset -\lambda$, by Ax5 we have $-\lambda \supset Y$, by CIC we infer $(Y \supset \sim X) \supset Y$ and, using Ax4 we get Y. By TD we conclude $(X \cap Y) \supset Y$.

Proof part c. Suppose $X\supset Y$, $X\supset Z$, $\sim Y\cup\sim Z$. By proposition 4 are derived $\sim Y\supset\sim X$, $\sim Z\supset\sim X$, applying proposition 5 is inferred $\sim X$, by TD results $(\sim Y\cup\sim Z)\supset\sim X$, by proposition 4 we can affirm $X\supset\sim(\sim Y\cup\sim Z)$, i.e., $X\supset\sim(\sim\sim Y\supset\sim Z)$, so $X\supset\sim(Y\supset\sim Z)$, and this means $X\supset\{Y\cap Z\}$. Applying TD 2 times concludes $(X\supset Y)\supset[(X\supset Z)\supset(X\supset\{Y\cap Z\})]$.
Proof part d. Suppose X, Y. Ax2 results $\lambda\supset X$, $\lambda\supset Y$, by part c derives $\lambda\supset(X\cap Y)$, using Ax1 infers $X\cap Y$. Applying TD 2 times concludes $X\supset [Y\supset(X\cap Y)]$.
Proof part e. Proposition 6 gives $(X\cap Y)\supset X$ and $(X\cap Y)\supset X$, using proposition 3 we derive $+[(X\cap Y)\supset X]$ and $+[(X\cap Y)\supset X]$, by proposition 2 we get $+(X\cap Y)\supset+X$ and $+(X\cap Y)\supset+X$, according to proposition 6 we conclude $+(X\cap Y)\supset(+X\cap+Y)$. To prove the reciprocal, by the part d we have $X\supset[Y\supset(X\cap Y)]\in FT$, using proposition 1 results $+\{X\supset[Y\supset(X\cap Y)]\}$, by proposition 2 we derive $+X\supset[+Y\supset+(X\cap Y)]$, again by proposition 2 and CIC we affirm $+X\supset[+Y\supset+(X\cap Y)]$, applying part c follows $(+X\cap+Y)\supset+(X\cap Y)$. Finally, applying part d and the definition of $\equiv$, we conclude $+(X\cap Y)\equiv(+X\cap+Y)$.

**Proposition 7.** The classical propositional calculus CPC with the language $\{\supset, \cap, \cup, \equiv, \sim\}$ is included in the propositional calculus GT.

Proof. Axioms 2, 3 and 4 along with propositions 4, 5 and 6, with the inference rule Mp determine CPC Rasiowa (1974).

**Proposition 8.** For $X,Y\in FT$. So, GT theorems: a) $--\lambda\supset\lambda$. b) $X\cup-X$. c) $\sim X\supset-X$. d) $-X\equiv\sim+X$. e) $+X\supset X$.

Proof part a. By Ax2 we have $\lambda\supset(--\lambda\supset\lambda)$, in addition to by Ax$\lambda$ of has $\lambda$, applying Mp we conclude that $--\lambda\supset\lambda$.
Proof part b. By Bx6 we have $\sim X\supset-X$, by definition it means $X\cup-X$.
Proof part c. By definition in Ax6.
Proof part d. By definition we have $\sim-X\equiv+X$, applying proposition 4 we conclude $-X\equiv\sim+X$.
Proof part e. By Bx6 we have $\sim X\supset-X$, applying CPC we deduce $\sim-X\supset X$, i.e., $+X\supset X$.

**Proposition 9.** For $X,Y\in FT$. So, GT theorems: a) $\sim+\sim X\equiv\otimes X$, $+\sim X\equiv\sim\otimes X$, $\sim+X\equiv-X$. b) $X\supset\otimes X$. c) $-X\supset\otimes\sim X$.
d) $\sim(Z_1\cap \dots \cap Z_k\cap Y)\in TT$ implies $\sim(+Z_1\cap\dots\cap+Z_k\cap \otimes Y)\in TT$.

Proof part a. By proposition 8 we have $-\sim X\equiv\sim+\sim X$, by definition it results $\otimes X\equiv\sim+\sim X$. By CPC we conclude $\sim\otimes X\equiv+\sim X$. By definition you have $\sim-X\equiv+X$, by CPC you get $-X\equiv\sim+X$,
Proof part b. By proposition 2 we have $+\sim X\supset\sim X$, using CPC we deduce $X\supset\sim+\sim X$, according to part a we conclude $X\supset\otimes X$.
Proof part c. By definition we have $-\sim X\supset\otimes\sim\sim X$, by CPC we conclude $-X\supset\otimes\sim X$.
Proof part d. Suppose $\sim(Z_1\cap \dots \cap Z_k\cap Y)\in TT$, which by CPC means, $(Z_1\cap\dots\cap Z_k)\supset\sim Y\in TT$. Using proposition 3 it turns out that $+((Z_1\cap\dots \cap Z_k)\supset\sim Y)\in TT$, from proposition 2 we infer $+(Z_1\cap \dots \cap Z_k)\supset+\sim Y\in TT$, by proposition 6 we get $(+Z_1\cap \dots \cap+Z_k)\supset+\sim Y\in TT$, which, by CPC implies $\sim(+Z_1\cap \dots \cap+Z_k\cap\sim+\sim Y)\in TT$, and for the part a, equivalent to $\sim(+Z_1\cap\dots\cap+Z_k\cap \otimes Y)\in TT$.

## III. Semantics

In this section, the semantics of possible worlds for the GT system are presented, in proposition 12, it is proved that the theorems of the GT system are valid formulas in the proposed semantics.

**Definition 4.** (S, Ma, <, V) is a model *for GT, it means that, S is a non-empty set* of possible worlds, Ma is a possible world, *called the* actual world, < is a binary relation in S, V is a valuation *of* $S\times FK$ at {0, 1}. The relationship**,** <, satisfies the following constraints. Reflexivity of <. RR: $(\forall M\in S)(M<M)$.

**Definition 5.** In the model Mo=(S, Ma, <, V), with $X,Y\in FT$.
V(M, X)=1 is abbreviated as M(X)=1, and means that in the possible world M, the formula X is *true*.
V(M, X)=0 is abbreviated as M(X)=0, and means that in the possible world M, the formula X is *false*.
X is true in Mo means that V(Ma, X)=1.

Valuation V satisfies the following rules: 1) V$\lambda$. M($\lambda$)=1. 2) V$\supset$. M($X\supset Y$)=1 equivalent to M(X)=1 implies M(Y)=1). 3) V$-$. M($-X$)=1 equivalent to $(\exists P\in S)(M<P$ and $P(X)=0)$.

**Proposition 10.** For $X,Y\in FT$. a) V$\sim$. M($\sim X$)=1 equivalent to M(X)=0.
b) V$\cup$. M($X\cup Y$)=1 equivalent to M(X)=1 or M(Y)=1.

c) V$\cap$. M($X\cap Y$)=1 equivalent to M(X)=M(Y)=1. d) V$\equiv$. M($X\equiv Y$)=1 equivalent to M(X)=M(Y).

e) $V+$. $M(+X)=1$ equivalent to $(\forall N\in S)(M<N$ implies $N(X)=1)$.

f) f) $V\otimes$. $M(\otimes X)=1$ equivalent to $(\exists P\in S)(M<P)$ $(P(X)=1)$.

g) $M(+X)=1$ implies $M(X)=1$. h) $M(-X)=0$ implies $M(X)=1$. i) $M(-\lambda)=0$

Proof parts a, b, c, and d. By CPC.

Proof part e. If $M(+X)=1$, then by part b, $V+$, follows $(\forall N\in S)(M<N$ implies $N(X)=1)$, in particular, $M<M$ implies $M(X)=1$, by the constraint RR we have $M<M$, therefore $M(X)=1$ .

Proof part f. If $M(\otimes X)=1$ equivalent to $M(\sim - X)=1$, by 1 equivalent to $M(-X)=0$, by $V\sim -$ equivalent to $(\forall P\in S)$ $(M<P$ implies $P(X)=1)$.

Proof part g. If $M(+X)=1$, for the part e, we have $(\forall N\in S)(M<N$ implies $N(X)=1)$, as $M<M$ then $M(X)=1$.

Proof part h. If $M(X)=0$, by $V-$ results $(\forall P\in S)(M<P$ implies $P(X)=1)$, but $M<M$, so $M(X)=1$.

Proof part i. If $M(-\lambda)=1$, by $V-$ affirms the existence of a world N, $M<N$ and $N(\lambda)=0$, which contradicts $V\lambda$.

**Definition 6.** For $X,X_1,..., X_n\in FT$, a formula X is said to be valid, denoted $X\in VT$, if and only if X is true in all models for GT, i.e., X is true in the current world of all models for GT. It is said that $\{X_1, ..., X_n\}$ *validates* Y if and only if $(X_1\cap X_2\cap ...\cap X_n)\supset Y\in VT$.

**Proposition 11.** For $X\in FT$. If X is an axiom of GT, then $X\in VT$.

Proof. Ax1. By $V\lambda$ we have for all $M\in S$, $V(\lambda)=1$. Hence, $Ax\lambda\in VT$.

Ax2, Ax3, Ax4. If X is one of the axioms Ax2, Ax3, Ax4, using the rule $V\supset$ and proceeding as usual for the validity of the intuitionistic propositional calculus in van Dalen (2004), it is concluded that $X\in VT$, i.e., Ax2, Ax3, $Ax14\in VT$.

Ax5. Suppose that $-\lambda\supset Z\notin VT$, so there is a model, such that in the present world M, $M(-\lambda\supset Z)=0$ by $V\supset$ results $M(-\lambda)=1$, using $V-$ it follows that there is $N\in S$, $M<N$ y $N(\lambda)=0$, which contradicts $V\lambda$. Hence, $Ax5\in VT$.

Ax6 Suppose that $(X\supset -\lambda)\supset -X\notin VT$, so there is a model, such that in the actual world M, $M((X\supset -\lambda)\supset -X)=0$, by $V\supset$ results $M(X\supset -\lambda)=1$ y $M(-X)=0$, applying $V-$ follows $(\forall P\in S)(M<P$ implies $P(X)=1)$, as $M<M$ follows $M(X)=1$, by $V\supset$ derives $M(-\lambda)=1$, using $V-$ it follows that there is $N\in S$, $M<N$ and $N(\lambda)=0$, which contradicts $V\lambda$. Hence, $Ax6\in VT$.

Ax7 Suppose that $[-(X\supset Y)\supset -\lambda]\supset[(-X\supset -\lambda)\supset(-Y\supset -\lambda)]\notin VT$, so there is a model, such that in the actual world M, $M([-(X\supset Y)\supset -\lambda]\supset[(-X\supset -\lambda)\supset(-Y\supset -\lambda)])=0$, by $V\supset$ results $M(-X\supset -\lambda)=1$ y $M(-Y\supset -\lambda)=0$, by $V\supset$ we get $M(-Y)=1$ y $M(-\lambda)=0$, applying $V\lambda$ we derive $M(-(X\supset Y))=0$. Since $M(-Y)=1$ according to $V-$ means $(\exists P\in S)(M<P$ y $P(Y)=0)$, since $M(-(X\supset Y))=0$, according to $V-$ we obtain $(\forall P\in S)(M<P$ implies $P(X\supset Y)=1)$, in addition, as $M(-X\supset -\lambda)=1$ y $M(-\lambda)=0$, by $V\supset$ we infer $P(-X)=0$, by $V-$ follows $(\forall P\in S)(M<P$ implies $P(X)=1)$, in particular $P(X)=1$, which is impossible. Therefore, $Ax7\in VT$.

Ax9 If $X\in\{Ax1, ..., Bx7\}$ then $-X\supset -\lambda$ is an axiom. Suppose that $X\in\{Ax1, ..., Ax7\}$ and $-X\supset -\lambda\notin VT$, so there is a model, such that in the present world M, $M(-X\supset -\lambda)=0$, by $V\supset$ results $M(-X)=1$, by V equivalent to $(\exists P\in S)(M<P$ y $P(X)=0)$, resulting in $X\notin VT$, but as $X\in\{Ax1, ... , Ax8\}$, the opposite has already been proved above. Hence, $Ax+\in VT$.

**Proposition 12.** For $X,Y\in FT$. a) $X\in TT$ then $X\in VT$. b) If $\{X1, ..., X_n\}$ implies Y then $\{X1, ..., X_n\}$ valid to Y.

Proof. Suppose $X\in T\Gamma$, $X\in VT$ is proved by induction over the length, L, of the proof of X.

Base step L = 1. It means that X is an axiom, which from proposition 11 follows that $X\in VT$.

Induction step. As an inductive hypothesis, we have that for every formula Y, if$\in$ $Y\in TT$ and the length of the proof of Y is less than L (where $L>1$) then $Y\in VT$. If $X\in TT$ and the length of the proof of X is L, then X is an axiom or X is a consequence of applying Mp in earlier steps of the proof. In the first case, we proceed as in the base case. In the second case, we have for some formula Y, proofs of Y and $Y\supset X$, where the length of both proofs is less than L, using the inductive hypothesis it is inferred that $Y\in VT$ and $Y\supset X\in VT$, so that, in the current world, M, of any model we have $M(Y)=1$ and $M(Y\supset X)=1$, by $V\supset$ it turns out that $M(X)=1$, consequently, $X\in V\Gamma$. Using the principle of mathematical induction, it has been proved that, for every $X\in FT$, $X\in TT$ implies $X\in VT$.

Suppose that $\{X_1, ..., X_n\}$ implies Y, applying CPC, we have $(X_1\cap ... \cap X_n)\supset Y\in TT$, from the part a is inferred, $(X_1\cap X_2\cap ...\cap X_n)\supset Y\in VT$, which means that $\{X_1, ..., X_n\}$ *validates* Y.

## IV. Semantic-deductive characterization

In this section, we present the characterization of GT with the semantics of the previous section. Completeness is proved in proposition 19 (valid formulas in semantics are theorems of GT), and characterization is achieved in proposition 20 (theorems of GT are the valid formulas of semantics and only they).

**Definition 7.** An *extension of a set of formulas C of d GT is obtained by altering the set of formulas* of C in such a way that the theorems of C are preserved, and that the language of the extension matches the language of GT. An extension is *locally consistent* if there is no $X \in FT$ such that both X and $\sim X$ are extension theorems. A set of formulas is *locally inconsistent* if a classical contradiction is derived from them, i.e., $Z \cap \sim Z$ is derived for some $Z \in FT$. An extension is *locally complete* if for all $X \in FT$, either X is an extension theorem or $\sim X$ is an extension theorem. To arrive at the demonstration of completeness in proposition 19, the strategy of the canonical model, presented in Henkin (1949), is followed.

**Proposition 13.** For $X \in FT$. a) GT is *locally* consistent. b) If $E \in EXT(GT)$, $X \notin TK\text{-}E$, and $E_x \in EXT(GT)$ are obtained by adding $\sim X$ as a new formula to E, then Ex is locally consistent.

Proof part a. Suppose GT not to be *locally* consistent, so there must be $Z \in FT$ such that $Z \cap \sim Z \in TT$, i.e. $Z \cap (Z \supset -\lambda) \in TT$, by CPC results $-\lambda \in TT$, by the validity theorem it is concluded that $-\lambda \in VT$, i.e., for $M \in S$, $M(-\lambda)=1$, this implies that there exists P, $M<P$ y $N(\lambda)=0$, which contradicts rule $V\lambda$. Therefore, GT is *locally* consistent.

Proof part b. For $X \notin TK\text{-}E$, and for Ex be the extension obtained by adding $\sim X$ as a new formula to E. Suppose that Ex is locally inconsistent, so that, for some $Z \in FT$, we have $Z, \sim Z \in TT\text{-}E_x$, by CPC comes $-\lambda \in TT\text{-}E_x$, by Ax5 $X \in TT\text{-}E_x$ is derived. But Ex differs from E only in that it has $\sim X$ as an additional axiom, so 'X is a theorem of Ex' is equivalent to 'X is a theorem of E from the set $\{\sim X\}$'. By TD it turns out that $X \supset \sim X \in TT\text{-}E$, and by CPC it is inferred that $X \in TT\text{-}E$, which is not the case, therefore $E_x$ is locally consistent.

**Proposition 14.** If $E \in EXT(GT)$ is locally consistent, then there exists $E' \in EXT(GT)$ which is locally consistent and complete.

Proof. For $X_0, X_1, X_2, \ldots$ an enumeration of all GT formulas. A sequence $E'_0, E'_1, E'_2, \ldots$ of extensions of E as follows: For $E'_0=E$. If $X_0 \in TT\text{-}E'_0$ is $E'_1=E'_0$, otherwise add $\sim X_0$ as a new formula to get $E'_1$ from $E'_0$. In general, given $t \geq 1$, to construct $E'_t$ from $E'_{t-1}$, proceed as follows: if $X_{t-1} \in TT\text{-}E'_{t-1}$, then $E'_t=E'_{t-1}$, otherwise for E't the extension of $E'_{t-1}$ obtained by adding $\sim X_{t-1}$ as a new formula. The proof is performed by mathematical induction over t.

Base step. t=0. Since E is consistent and $E'_0=E$, by hypothesis, it follows that $E'_0$ is locally consistent.

Inductive step. Inductive hypothesis: $E'_{t-1}$ is locally consistent with $t \geq 1$. By proposition 4, $E'_t$ is locally consistent. By the principle of mathematical induction, it follows that all $E'_t$ is locally consistent.

E' is defined as that extension of E, which has as new formulas those formulas that are new formulas of at least one of the $E'_t$. If E' is not *locally* consistent, then there exists $X \in FT$ such that, $X, \sim X \in TT\text{-}E'$, but the proofs of X and $\sim X$ in E' are finite sequences of formulas, so that each proof can only contain particular cases of a finite number of axioms or new formulas of E', so there must be a *t*, large enough that all these new formulas used are elements of $E'_t$, resulting in $X, \sim X \in TT\text{-}E'_t$, which is impossible since $E'_t$ is locally consistent. Therefore, E' is locally consistent.

To prove that E' is complete, for $X \in FT$. X must appear in the list $X_0, X_1, X_2, \ldots$, suppose X is $X_k$. If $X_k \in TT\text{-}E'_k$, then $X_k \in TT\text{-}E'$, since $E' \in EXT(E'_k)$, if $X_k \notin TT\text{-}E'_k$, then according to the construction of $E'_{k+1}$, $\sim X_k$ is a new formula of $E'_{k+1}$, so $\sim X_k \in TT\text{-}E'_{k+1}$, and then $\sim X_k \in TT\text{-}E'$. Thus, in any case, we have $X_k \in TT\text{-}E'$ or $\sim X_k \in TT\text{-}E'$, so E' is locally complete.

**Proposition 15**. For $Y, Z_1, \ldots, Z_k \in FT$. If $\{+Z_1, \ldots, +Z_k, Y\}$ is locally consistent then $\{Z_1, \ldots, Z_k, \otimes Y\}$ is locally consistent.

Proof. Suppose $\{Z_1, \ldots, Z_k, Y\}$ is locally inconsistent in GT, so there exists a formula $W \in FT$ such that, from $\{Z_1, \ldots, Z_k, Y\}$, $W \cap \sim W$ is inferred in GT, using CPC it turns out that $\sim(Z_1 \cap \ldots \cap Z_k \cap Y) \in T\Gamma$, by fact 1 is derived $\sim(+(Z_1 \cap \ldots \cap Z_k) \cap \otimes Y) \in TT$, by proposition 9 $\sim(+Z_1 \cap \ldots \cap +Z_k \cap \otimes Y) \in TT$ so $\{+Z_1, \ldots, +Z_k, \otimes Y\}$ is locally inconsistent in GT. It has been proved that $\{Z_1, \ldots, Z_k, Y\}$ locally inconsistent implies that $\{+Z_1, \ldots, +Z_k, \otimes Y\}$ locally inconsistent, i.e., $\{+Z_1, \ldots, +Z_k, \otimes Y\}$ locally consistent implies $\{Z_1, \ldots, Z_k, Y\}$ locally consistent.

**Definition 8.** Be locally consistent and complete E,F∈EXT(GT). F is said to be subordinate to E if and only if there is Y∈FT, such that Y∈E, and furthermore for every ⊗Z∈FT, such that +Z∈E, we have to Y,Z∈F.

**Proposition 16.** For E∈EXT(GT), X∈FT. If E is locally consistent and complete and ⊗X∈E, then there exists F∈EXT(GT) locally consistent and complete such that X∈F and F is subordinate to E.

Proof. Suppose ⊗X∈E. For $E_X$={X}∪{Z: +Z∈E}, since E is locally consistent, then for proposition 15, $E_X$ is also locally consistent. By adding to $E_X$ the axioms of GT and all their consequences, we get an extension of GT that includes $E_X$, using proposition 14, we construct a locally consistent and locally complete extension F of GT which includes $E_X$. Like X∈$E_X$, also X∈F. If +W∈E, by definition, W∈$E_X$, so W∈F. Therefore, F is subordinate to E.

**Proposition 17**. For locally consistent and complete *E,F,G*∈EXT(GT). RR. Reflexivity. F is subordinate to F.

Proof. For X be the axiom Ax1, so X∈TT, and since by proposition 9 we have X⊃⊗X, it follows that ⊗X∈TT, then X∈F and ⊗X∈F. Suppose that +W∈F, by proposition 8 ∈it follows that W∈F. Hence, F subordinate to F.

**Proposition 18.** If E'∈EXT(GT) is *locally* consistent, then there exists a model in which all X∈TT-E' is true.

Proof. The model (S, ME, <, V) is defined as follows: For E, F, G, ..., be *locally* consistent and complete extensions of E' (E the initial and the other subordinates), presented in the preceding propositions. To each extension F, a possible world MF is associated, for S the set of such possible worlds and ME the actual world. The accessibility relation, <, is constructed as follows: MF<MG if and only if G is subordinate to F.

For each MF∈S and for each X∈F, V(MF,X)=1 if X∈F and V(MF,X)=0 if ~X∈F, where F is the locally consistent and complete extension associated with MF. Note that V is functional because F is *locally* consistent and complete.

To claim that M is a model, rules 1 to 3 of definition 5 must be guaranteed.

1. By Axλ You have λ∈TT, so λ∈F, i.e., V(MF,λ)=1. Therefore, Vλ is satisfied.

2. In the case of the conditional X⊃Y. Using CPC we have the following chain of equivalences: V(MF,X⊃Y)=0, i.e. ~(X⊃Y)∈F, by CPC we follow X∩~Y∈F, resulting in CPC that X∈F y ~Y∈F, which means that V(MF,X)=1 and V(MF,Y)=0, so V⊃ is satisfied.

3. In the case of rule V−. Sea MF is a world associated with F, MG is a world associated with G and Z∈FK. Suppose that V(MF,−Z)=1, so −Z∈F, and by proposition 9 ⊗~Z∈F, by proposition 16, there exists G subordinate to F, such that ~Z∈G, resulting that, (∃MG∈S)(MF<MG y MG(Z)=0).

To prove the reciprocal, suppose (∃MG∈S)(MF<MG y MG(Z)=0). If V(MF,−Z)=0, then it follows that ~−Z∈F, i.e., +Z∈F, and since MF<MG, i.e., G is subordinate to F, then Z∈G, i.e., MG(Z)=1, result, by the hypothesis, that G is locally inconsistent, which is not the case. Therefore, V(MF, −Z)=1. Since the reciprocal has already been proved, then definition V− is satisfied.

Based on the above analysis, it is inferred that V is a valuation, and since the constraint RR is guaranteed by proposition 17, it is finally concluded that M is a model.

To conclude the proof, for X a theorem of E', so X is in E'. Therefore, using the definition of V, it turns out that V(ME,X)=1, i.e., X is true in the model M=(S, ME, <, V).

**Proposition 19.** For X,$X_1$, ..., $X_n$∈FT. a) If X∈VT then X∈TT.

b) If {$X_1$, …, $X_n$} validates Y then Y is a consequence of {$X_1$, …, $X_n$}.

Proof part a. If X∉TT, then, by proposition 13, the extension E', obtained by adding ~X as a new formula, is locally consistent. Thus, according to proposition 18, there is a model M such that every theorem of E' is true in M, and since ~X∈TT-E', then ~X is true in M, i.e., X is false in M, hence X∉VT. It has been proved that X∉TT implies X∉VT, i.e., X∈VT implies X∈TT.

Proof part b. Suppose {$X_1$, …, $X_n$} validates Y, i.e., ($X_1$∩$X_2$∩…∩$X_n$)⊃Y∈VT, by part a, follows that, ($X_1$∩$X_2$∩…∩$X_n$)⊃Y∈TT. If {$X_1$, …, $X_n$} are assumed, by CPC Y is inferred, therefore Y is a consequence of {$X_1$, ..., $X_n$}.

**Proposition 20.** For X,Y,$X_1$,..., $X_n$ ∈FT. a) X∈VT if and only if X∈TT.

b) {$X_1$, ..., $X_n$} validates Y if and only if {$X_1$, ..., $X_n$} implies Y.

Proof. Consequence of propositions 12 and 19.

## V. Existential graphs

In this section, we present the primitive existential graphs for the GT system. For the construction of existential graphs, a variant of the notation proposed by Peirce in 4.378 of Peirce'*s Collected Papers* (1965) is used.

**Definition 9.** The set, GET, of existential graphs for the GT system, is constructed from a set of atomic graphs, GA, and the constant (λ empty graph, ='_'), as follows.
P∈GA implies P∈GET. λ∈GET. X∈GET implies {X}∈GET. X,Y∈GET implies (X(Y))∈GET.

**Definition 10.** On the graph (X(Y)) it is called *a conditional graph.* The part (…) is called *the outer cut of the conditional*. The other part is called, it's called *the internal cut of the conditional*. In (X(Y))), X is called antecedent and Y *consequent*. Conditional cuts are called *continuous cuts.*
In the {Z} graph, the *{...}*, it's called *a broken cut.*

**Definition 11.** For X,Y,Z∈GET. A graph X is said to be in an even *region*, denoted $X_p$, if X is surrounded by an even number of *cuts (continuous* and/or *broken).* X is in an odd region, denoted $X_i$, if X is surrounded by an odd number of *cuts (continuous* and/or *broken).*
$X_{nc}$ means that graph X is in a region surrounded by *n* continuous and/or broken slices (n=*0, 1, 2, 3, ...), where* n *can be* odd or even.
$X_{ncc}$ means that X is in a *region of continuous cuts only*, i.e., no broken cuts appear.
$X_{1cq}$ means that X is in a region with at least one broken cut.
*X stands for ({X}), and *X is said to be a *strong graph*.

**Definition 12.** For X∈GET. Lambda is defined as the assertion sheet λ = '_'. Strong statement is defined as *X = ({X}) Total falsehood is defined as = {⊥}.

**Definition 13.** RTRA Primitive Transformation Rules
R1. Strong double cut writing. The strong double cut is a graphical theorem. ({_}), ({λ}).
R2. Graphics erasure. A chart can be *deleted* when it is in *an even* region. $XY_p|\Rightarrow X_p$
Writing charts. In an *odd* region, any chart *can be* written. $X_i|\Rightarrow XY_i$
R3. Unrestricted iteration and de-iteration of graphics in continuous slice-only region. A *chart* can *be iterated* or *de-iterated* in any region, odd or even, if the region is only surrounded by zero or more continuous slices.)
$Y(X)_{ncc} \Leftrightarrow Y(XY)_{ncc.}$ X ⇔ XX $\quad$ $*Y(X)_{1cq} \Leftrightarrow *Y(X*Y)_{1cq.}$
R4. Erased in a cut. A continuous cut can be *partially erased* (generating a broken cut) when it is in an *even* region. $(X)_{p|}\Rightarrow \{X\}_p$
Writing on a cut. A broken cut can be *completed* (generating a continuous cut) when it is in an *odd* region.
$\{X\}_i|\Rightarrow (X)_i$
R6. Erasure of continuous double cutting. A continuous double cut can be *erased* in *an even* region. $((X_p))|\Rightarrow X_p$
Continuous double cut writing. A continuous double cut can be *written around* a graph that is in *an odd* region.
$X_i|\Rightarrow((X_i))$

Rules 7, 8 and 9 are called implicit rules, since, given their obviousness and graphic naturalness, they may not be referenced, but they are applied.
R7. Concatenation. Two graphs that are in the same region can be concatenated. Conversely, two graphs that are concatenated can be separated in the same region. X, Y ⇔ YX, in any region.
R8. Commutativity. Two concatenated charts can be rewritten by changing the order. XY ⇔ YX, in any region.
R9. Associativity. In three graphs that are concatenated, the order in which they were concatenated is irrelevant. Initially, the first is concatenated with the second and this result is concatenated with the third, or the first is concatenated with the result of concatenating the second with the third. XY, Z ⇔ X, YZ XYZ, ⇔ in any region.

**Definition 14.** For X∈GET. X is a *graphical theorem* of GET, denoted X∈TGET, if there is *a proof* of X from the graph, using the graph λ transformation rules, i.e., X is the last row of a finite sequence of lines, in which each of the lines is, or is inferred from previous rows, using the transformation rules. Or to put it briefly, X∈TGET if and only if λ>>X. The number of lines, of the finite sequence, is referenced as the *length* of the proof of X.
Y>>X, means that X is obtained from Y using a finite number of transformation rules

**Proposition 21.** For X,Y∈GET. $(\forall R\in RTRA)[\, X_p\overset{R}{\Rightarrow}Y](\exists R'\in RTRA)[\, Y_i\overset{R'}{\Rightarrow}X]$

Proof: Simple inspection of the primitive rules.

**Proposition 22.** For X,Z∈GET. a) $X>>Z \mid\!\Rightarrow X_p>>Z$. b) $X>>Z \mid\!\Rightarrow Z_i>>X$.

Proof: For X,Z∈GET, suppose X>>Z. it must be proved that $[X_p>>Z$ and $Z_i>>X]$.

If X>>Z then there are $R_1, \ldots, R_n\in RTRA$, and there are $X_1, \ldots, X_{n-1}\in GET$, such that $XR_1X_1R_2X_2\ldots X_{n-1}R_nZ$, and the *length* of the X>>Z transformation is said to be n and $X>>_nZ$ is denoted. The proof is performed by induction on the length of the transformation.

Base step. n=1. It means that only one of the primitive rules was applied, and since X is in an even region, then R must be of the form $X_p\overset{R}{\Rightarrow} Z$ with R∈RTRA. From proposition 21 it is inferred that $Z_1\overset{R'}{\Rightarrow} X$ with R'∈RTRA.

Inductive step. Inductive hypothesis $(\forall n>1)[W>>_nK \Rightarrow\{W_p>>K$, so $X>>_nX_n$ and $X_nR_{n+1}Z$. Applying the inductive hypothesis and proposition 21 we get $X_p>>X_n$ and $X_nR_{n+1}Z$, $X_{ni}>>X$ and $Z_iR'_{n+1}X_n$. So, $X_p>>Z$ y $Z_i>>X$.

By the principle of mathematical induction, the truth of the proposition is concluded.

**Proposition 23.** For X,Z∈GET. $X>>Z \mid\!\Rightarrow(X(Z))$.

Proof. Suppose X>>Y. $\lambda \overset{R1}{\Rightarrow} (\{_\}) \overset{R4}{\Rightarrow} ((_)) \overset{R2}{\Rightarrow} (X(_)) \overset{R4}{\Rightarrow} (X(X)) \overset{X\gg Y \text{ and proposition 22}}{\Longrightarrow} (X(Y))$. Hence, $X>>Y \Rightarrow (X(Y))$.

**Proposition 24.** For X∈GET. a) $((_))$, $*\lambda$, $\{(_)\}$, $\otimes\lambda$, $\lambda$. b) $X_p \mid\!\Rightarrow((X_p))$. c) $((X_i)) \mid\!\Rightarrow X_i$. d) $((X)) \Leftrightarrow X$.

Proof part a. By R1 we have $(\{_\})$, i.e. $(\{\lambda\})$, which means $*\lambda$. By having $(\{_\})$, by R4 we infer $((_))$, again by R4 we derive $\{(_)\}$, i.e. $\{(\lambda)\}$, which means $\otimes\lambda$. Since we already have $((_))$, i.e. $((\lambda))$, Using R6 we conclude $\lambda$.

Proof part b. Suppose X, for part a, we have $((_))$, applying R3 we conclude ((X)). X⇒((X)) has been tested. Part c is obtained by proposition 22. Part d results from parts b and c together with rule R6.

**Proposition 25.** For X∈GET. a) $*X_p \mid\!\Rightarrow X_p$. $(\{X_p\}) \mid\!\Rightarrow X_p$. b) $X_i \mid\!\Rightarrow *X_i$, $X_i \mid\!\Rightarrow(\{X_i\})$.

Proof. Suppose $(\{X_p\})$, by R4 we get $((X_p))$, according to R6 we derive $X_p$. Therefore, $(\{X_p\}) \mid\!\Rightarrow X_p$. Using proposition 22 concludes $X_i \mid\!\Rightarrow (\{X_i\})$.

**Proposition 26.** For X∈GET. $\{_\}_p \mid\!\Rightarrow X$, for everything X∈GET

Proof: By R1 we have $(\{_\})$, using R2 we infer $(\{_\}(X))$, Suppose $\{_\}$, by R3 we ensure ((X)), according to R6 we conclude X, we have proved, $\{_\}>>X$. Proposition 22 concludes $\{_\}_p \mid\!\Rightarrow X$.

**Proposition 27.** For X∈GET. a) X∈TGET implies $X_p \mid\!\Rightarrow *X_p$. b) X∈TGET implies $X_p \mid\!\Rightarrow(\{X\})$.

Proof. If X∈TGET then $\lambda>>X$, applying proposition 22 we derive $\lambda_p>>X_p$, and since $\lambda$ by proposition 12 we have $(\{\lambda_p\})$, then we conclude $(\{X_p\})$. Therefore, $X_p \Rightarrow (\{X\})$, i.e. $X_p \Rightarrow *X$.

**Proposition 28.** For X∈GET. a) $X>>\{_\}$ implies **{X}**. b) **{X}>>{_}** implies X.

Proof. Part a. Suppose $X>>\{_\}$, by TD results $(X(\{_\}))$, by proposition 24 we have $(\{_\})$, applying R3 we deduce (X), according to R4 we conclude {X}.

Part b. Suppose $\{X\}>>\{_\}$, by TD results $(\{X\}(\{_\}))$, by proposition 24 we have $(\{_\})$, applying R3 we deduce ({X}), according to proposition 25 we conclude X.

## VI. Equivalence between GT and GET

In this section, the equivalence between GT and GET is presented, initially, in proposition 32, it is proved that the theorems of GT are graphical theorems of GET, in proposition 37, it is proved that the graphical theorems of GET are valid in the semantics of possible worlds, in proposition 40, it is proved that the theorems of GT are exactly the graphic theorems.

**Definition 15.** FT translation function [_]' in GET. Be X,Y∈FT and P∈FA. 1) P' = P. 2) $[X\supset Y]' = (X'(Y))$. 3) $[X\cup Y]' = ((X')(Y'))$. 4) $[-X]' = \{X'\}$. 5) $[X\cap Y]' = X'Y'$. 6) $[\sim X]' = (X')$. 7) $\lambda' = \lambda$.

**Proposition 29.** For X∈FT. If X is an axiom of GT then X'∈TGET.

Proof. Using primitive rules, you have:

Ax1. **$\lambda$**. Since $\lambda'=\lambda$, then for proposition 24, $\lambda$ is also a graphical theorem.
Ax2. $X\supset(Z\supset X)$. By proposition 24 we have ((_)) and ((_)), by R3 we derive ((((_)))), using the rule R2 follows (X'(((_))), applying R3 results (X'(((X')))), again by R2 we derive (X'((Z'(X')))), i.e. (Ax2)' is a graphical theorem.
Ax3. $(X\supset(Y\supset Z))\supset((X\supset Y)\supset(X\supset Z))$. Suppose (X'((Y'(Z')))), (X'(Y')) and X'. By R3 we can deduce ((Y'(Z')))), ((Y')), applying R6 we follow (Y'(Z')), Y', using R3 we infer ((Z')), by R6 we conclude Z'. By TDG 2 times, it has been proved that (Ax3)' is a graphical theorem.
Ax4. $[(X\supset Y)\supset X]\supset X$. Suppose ((X'(Y'))(X')), by R2 we get ((X')(X')), as R3 results ((X')), applying R2 derives X'. Using TDG it is concluded that (Ax4)' is a graphical theorem.
Ax5. $-\lambda\supset Z$. By proposition 26 Tri rule we have {_}>>Z, applying TDG we conclude that (Ax5)' is a graphical theorem.
Ax6. $(X\supset-\lambda)\supset-X$. Suppose (X'({_})), whence X'>>{_}, by proposition 30 is derived {X'}. Using TDG it is concluded that (Ax6)' is a graphical theorem.
Ax7. $[-(X\supset Y)\supset-\lambda]\supset[(-X\supset-\lambda)\supset(-Y\supset-\lambda)]$. Suppose ({(X'(Y'))}({_})), ({X'}({_})). i.e., *(X'(Y')), *X', by proposition 25 follows *(*X'(Y')), using R5 we deduce *((Y')), applying R2 we affirm *Y'. Using TDG it is concluded that (Ax7)' is a graphical theorem.
Ax+. If $X\in\{Ax1, \ldots, Ax7\}$ then $-X\supset-\lambda$ is an axiom. Consequence of proposition 27: $X\in TGET$ implies $X_p \Rightarrow {}^*X_p$.

**Proposition 30.** For $X\in FT$. a) If $X\in TT$ then $X'\in TGET$. b) If X⟩⟩And then X'>>Y'.
Proof part a. Induction on the length of the X demonstration in GT.
Base step. If the length of the proof is 1, then X is an axiom, by proposition 29 $X'\in TGET$.
Induction step. The inductive hypothesis is: if $Y\in TT$ and the length of the proof of Y is less than L then $Y'\in TGET$. Suppose $X\in TT$ and that the length of the proof of X is L, so X is an axiom or obtained from previous steps using Mp. In the first case, proceed as in the base step. In the second case, Y and $Y\supset Z$ are taken in previous steps of the proof of X, i.e., the lengths of the proofs of Y and $Y\supset X$ are less than L, by the inductive hypothesis it turns out that $Y'\in TGET$ and $(Y'(X'))\in TGET$, applying R3 infers $((X'))\in TGET$, using R6 we conclude $X'\in TGET$.
By the principle of mathematical induction, it is proved that the theorems of GT are graphic theorems.
Proof part b. If X implies Y, then $X\supset Y\in TT$, by the part a, $(X'(Y'))\in TGET$, i.e., $\lambda$>>(X'(Y')), if X' is assumed, by R3 follows ((Y')), applying R6 results in Y', so X'>>Y'.

**Definition 16.** Translation function, (_)'' of GET in FT. For $X,Y\in GET$, $P\in GA$. 1) P''=P. 2) $\lambda''=\lambda$.
3) $(X(Y))''=X''\supset Y''=[X''\supset[Y'']]$. 4) $((X)(Y))''=X''\cup Y''=[[X'']\cup[Y'']]$. 5) $\{X\}''=-X''=-[X]''$. 6) $(XY)''=X''\cap Y''$.
Observation. The notation [X''] simply means the formula X'', indicates that the associated graph X is surrounded by 1 slice. Therefore, the number of slices (_) surrounding a chart in GET matches the number of square brackets [_] surrounding the associated formula in GT.

**Proposition 31.** For $X,Y,Z\in GT$. If $X\Rightarrow Y$ then Y'' is validly inferred from X'', in the case of the most elementary versions of the rules. Generalization to arbitrary odd or even regions will be presented later.
Proof. R1. ({_}). If $\sim-\lambda$ it is invalid then there is a model with the actual world M such that, $M(\sim-\lambda)=0$, by V∼ follows $M(-\lambda)=1$, by V− it follows that $(\exists P\in S)(M<P \text{ y } P(\lambda)=0)$. Which contradicts V$\lambda$. Therefore, R1'' is valid.
R2a. XZ|⇒ X. Consider an arbitrary model with the actual world M. Suppose that $M(X''\cap Y'')=1$, by V∩ we derive M(X')=1. Therefore, R2'' is a valid rule in GT.
R3. Y(X)⇔Y(XY). Consider an arbitrary model with the actual world M. Suppose that $M(Y''\cap\sim X'')=1$, by V∩ follow M(Y'')=1, by V follow M(Y'')=1, $M(\sim X'')=1$, according to V∼ derive M(X'')=0, using V∩ we affirm $M(X''\cap Y'')=0$, by V∩ we infer $M(\sim(X''\cap Y''))=1$, Applying V∼ we get $M(Y''\cap \sim(X''\cap Y''))=1$. It has been proved that $M(Y''\cap\sim X'')=1$ implies $M(Y''\cap\sim(X''\cap Y''))=1$.
To prove the reciprocal, suppose that $M(Y''\cap\sim(X''\cap Y''))=1$, by V follow M(Y'')=1 and $M(\sim(X''\cap Y''))=1$, resulting $M(X''\cap Y'')=0$, i.e. M(Y'')=0 or M(X'')=0, but M(Y'')=1, so M(X'')=0, i.e. $M(\sim X'')=1$, and by V∩ we deduce $M(Y''\cap\sim X'')=1$. It has been proved that $M(Y''\cap\sim(X''\cap Y''))=1$ implies $M(Y''\cap\sim X'')=1$. Therefore, R3'' is a valid rule in GT.
R4a. (X)|⇒{X}. Consider an arbitrary model with the actual world M. Suppose that $M(\sim X'')=1$, i.e. M(X'')=0, and since $M<M$, we can affirm $(\exists N)(M>N)(N(X'')=0$, by V− we derive $M(-X'')=1$. It has been proved that $M(\sim X'')=1$ implies $M(-X'')=1$. Therefore, R4'' is a valid rule in GT.

R6a. ((X))|⇒ X. Suppose that M(~~X")=1 is equivalent to M(~X")=0, again by the same rule we conclude M(X")=1. It has been proved that M(~~X")=1 implies M(X")=1. Therefore, R6" is a valid rule in GT.

**Proposition 32.** For X,Y,Z,W,V∈FT. If X⊃Y∈VT are valid in GT: a) (X∩Z)⊃(Y∩Z). b) −(Y∩Z)⊃−(X∩Z). c) −(W∩−(X∩Z))⊃−(W∩−(Y∩Z)).

Proof. Suppose X ⊃Y∈VT, so you have the initial result, for every model, (S, Ma, <, V), and for every M∈S, M(X⊃Y)=1. Proof part a. Is CPC result.

Proof part b. Be M be the actual world of any model. Suppose that M(−(Y∩Z))=1, then by V−, there is P∈S, M<P and P(Y∩Z)=0. Suppose M(−(X∩Z))=0, times V− it follows that, for every N∈S, M<N implies N(X∩Z)=1, and as M<P then P(X∩Z)=1, times V∩ we have P(X)=1 and P(Z)=1, as P(Y∩Z)=0, by V∩**,** P(Y)=0 is derived, but as P(X)=1, and by the initial result, P(X⊃Y)=1, then by V⊃, we get P(Y)=1, which is not the case, so M(−(X∩Z))=1. It has been proven that, M(−(Y∩Z))=1 implies M(−(X∩Z))=1, which by V⊃ means that M(−(Y∩Z)⊃−(X∩Z))=1. Therefore, −(Y∩Z)⊃−(X∩Z)∈VT.

Proof part c. Be M be the actual world of any model. Suppose that M(−(W∩−(X∩Z)))=1, then by V−, there is P∈S, M<P and P(W∩−(X∩Z))=0. Assumption 2, M(−(W∩−(Y∩Z)))=0, times V− it follows that, for every N∈S, M<N implies N(W∩−(Y∩Z))=1, and as M<P then P(W∩−(Y∩Z))=1, by V∩ we have P(W)=1 and P(−(Y∩Z))=1, by V− it follows that, there is Q∈S, P<Q and Q(Y∩Z)=0, such as P(W)=1 and P(W∩−(X∩Z))=0, times V∩ P(−(X∩Z))=0, which is by V− it follows that, for each N∈S, P<N implies that N(X∩Z)=1, as P<Q, in particular, Q(X∩Z)=1, times V∩ followed by Q(X)=1 and Q(Z)=1, for the initial result, we have Q(X⊃Y)=1, times V⊃ Q(Y)=1 follows, and Q(Z)=1 is followed by V∩ Q(Y∩Z)=1 is derived, which is not the case, so M(−(W∩−(Y∩Z)))=1. It has been proven that, M(−(W∩ −(X∩Z)))=1 implies M(−(W∩−(Y∩Z)))=1, which by V⊃ means that M(−(W∩−(X∩Z))⊃−(W∩−(Y∩Z)))=1. Therefore, −(W∩−(X∩Z))⊃−(W∩−(Y∩Z))∈VT.

**Proposition 33**. for X,Y,Z,W,V∈FΓ. If X⊃Y∈TT are GT theorems: a) (X∩Z)⊃(Y∩Z). b) −(Y∩Z)⊃−(X∩Z). c) −(W∩−(X∩Z))⊃−(W∩−(Y∩Z)).

Proof. A direct consequence of propositions 20 and 32.

**Proposition 34.** For X,Y,Z,W,V∈FT. If X⟩⟩Y then, a) $X_p$⟩⟩Y. b) $Y_i$⟩⟩X.

Proof part a. Induction in the number, *n*, of negations surrounding X.

Base step. n=0. X∩Z⟩⟩Y∩Z, is satisfied by proposition 21.

n=1. There are 2 possibilities, –{X∩Z} and ~{X∩Z}. For proposition 33, we have –{Y∩Z}⊃–{X∩Z}, for CPC we have ~{Y∩Z}⊃~{X∩Z}, so the proposition is satisfied when n=1.

n=2. There are 4 possibilities, –{W∩–{X∩Z}}, –{W∩~{X∩Z}}, ~{W∩–{X∩Z}} y ~{W∩~{X∩Z}}. By proposition 33, we have that, –{W∩–{X∩Z}}⊃–{W–∩{Y∩Z}} y –{W∩~{X∩Z}}⊃–{W∩~{Y∩Z}}, the other 2 cases ~{W∩–{X∩Z}}⊃~{W–∩{Y∩Z}} y ~{W∩~{X∩Z}}⊃~{W∩~{Y∩Z}} are taken for CPC, so the proposition is satisfied when n=2.

Inductive step. Part a. As an inductive hypothesis we have that, if X is surrounded by 2n negations, then X>>Y. By proposition 33 we have that, –{W∩–{X∩Z}}⊃–{W∩–{Y∩Z}} y –{W∩~{X∩Z}}⊃–{W∩~{Y∩Z}}, and by CPC we have ~{W∩–{X∩Z}}⊃~{W–∩{Y∩Z}} y ~{W∩~{X∩Z}}⊃ ~{W∩~{Y∩Z}}, in the region surrounded by 2n negations, and these are the only cases for which two other negations can be added to X. Therefore, if X is surrounded by 2n+2 slices, i.e., by 2(n+1) slices, then X>>Y.

Proof part b. As an inductive hypothesis we have that, if X is surrounded by 2n+1 negations, then Y>>X. By proposition 33 we have that, –{W∩–{X∩Z}}⊃–{W∩–{Y∩Z}} y –{W∩~{X∩Z}}⊃–{W∩~{Y∩Z}}, and by CPC we have ~{W∩–{X∩Z}}⊃~{W–∩{Y∩Z}} y ~{W∩ ~{X∩Z}}⊃~{W∩~{Y∩Z}}, in the region surrounded by 2n+1 negations, and these are the only cases for which two other negations can be added to X. Therefore, if X is surrounded by 2n+1+2 negations, i.e. by 2(n+1)+1 negations, then Y>>X.

By the principle of mathematical induction, the proposition has been proved.

**Proposition 35.** For X,Y,Z∈GET. a) Primitive GET rules are valid rules in GT semantics. b) If X∈GET then X"∈VG.

Proof part a. Direct consequence of propositions 31 and 34.

Proof part b. If X∈TT then λ>>X, then there are R1, ..., Rn RTRA, and there are $X_1$, ..., $X_{n-1}$∈GET, such that ∈λ$R_1X_1R_2X_2$... $X_{n-1}R_nX$.

The proof is performed by induction over the length L of the demonstration.

PB. Base step. L=1. It means that only one of the primitive rules was applied, then X"∈VG.

PI. Inductive step. Inductive hypothesis: The proposition is valid if L<n+1 with n>0. Be L=n+1, so $\lambda R_1X_1R_2X_2 \ldots X_{n-1}R_nX_nR_{n+1}X$, i.e., $\lambda R_1X_1R_2X_2 \ldots X_{n-1}R_nX_n$ and $X_nR_{n+1}X$, both demonstrations with a length shorter than n+1. Applying the inductive hypothesis, it turns out that $X_n$''∈VG and from $X_n$'' is validly inferred X", *hence X"*∈VG.
By the principle of mathematical induction, the truth of the proposition is concluded.

**Proposition 36.** For X,Y∈GET. a) If X∈TGET then X"∈TT. b) If X>>Y then X" implies Y".
Proof part a. By proposition 20 we have that, X"∈VT if and only if X"∈TT, and by proposition 35 we have that, if X∈TGET then X"∈VG. Therefore, if X∈TGET then X"∈TT.
Proof part b. Consequence of part a and proposition 20.

**Definition 17.** Be T1=[_]':FG→GG and T2=[_]":GG→FG, be the translation functions presented in definitions 15 and 16. They are defined: the composite function T1oT2:GG→GG such that (T1oT2)[X]=T1[T2[X]], the composite function, T2oT1:FG→FG such that (T2oT1)[X]=T2[T1[X]], the identity function in FG, $Id_{FK}$:FG→FG such that ($Id_{FG}$)[X]=X, the identity function in G $Id_{GG}$:GG→GG such that ($Id_{GG}$)[X]=X.

**Definition 18.** For P∈GA, X,Y∈GET. The function, C*, complexity of* a graph, assigns each graph a non-negative integer, as follows: 1) C[P] = C[λ] = 0. 2) C[{X}] = 1+C[X]. 3) C[XY] = 1+max{C[X], C[Y]}.
4) C[[(X)(Y)]] = 2+max{C[X], C[Y]}. 5) C[(X(Y))] = 1+max{C[X], C[Y]+1}.

**Definition 19.** Sean P∈FA; X,Y∈FT. The function, K*, complexity of a formula,* assigns each formula a non-negative integer, as follows: 1) K[P] = K[λ] = 0. 2) K[−X] = 1+K[X].
3) K[X∩Y] = K[X∪Y] = K[X⊃Y] = 1+max{K[X], K[Y]}.

**Proposition 37.** For P∈GA, P∈FA, G,G1,G2∈GET and X,X1,X2∈FT. a) T1oT2 = $Id_{GK}$. b) T2oT1= $Id_{FG}$.
c) T1 is the inverse function of T2. d) T2 is the inverse function of T1.
Proof part a. Induction on the complexity, C, of graph G.
Base step. C[G]=0, then there are 2 cases.
Case 1: G=P. (T1oT2)[P]=T1[T2[P]]=T1[P]=P. Case 2: G=λ. (S1oS2)[λ]=T1[T2[λ]]=T1[λ]=λ.
Inductive step. C[G]≥1. As an inductive hypothesis we have that (T1oT2)[G1]=G1, (T1oT2)[G2]=G2.
There are 4 cases. Case 3: G={G1}. (T1oT2)[{G1}]=T1[T2[{G1}]]=T1[T2[G1]]={T1[T2[G1]]}={−G1}.
Case 4: G=G1G2. (T1oT2)[G1G2]=T1[T2[G1G2]]=T1[T2[G1]T2[G2]]=T1[T2[G1]]∩T1[T2[G2]]=G1G2.
Case 5: G=(G1(G2)). (T1oT2) [(G1(G2))]=T1[T2[(G1(G2))]]=T1(T2[G1] ⊃T2[G2])
=(T1[T2[G1]](T1[T2[G2]]))=(G1(G2)).
Case 6: G=((G1)(G2)). (T1oT2) [((G1)(G2))]=T1[T2[((G1)(G2)]]=T1[T2(G1)T2(G2)]
=((T1[T2[G1])(T1[T2[G2]))=((G1)(G2)).
By the principle of mathematical induction it has been proved that (T1oT2)=$Id_{GG.}$
Proof part b. Induction on the complexity, K, of the formula X.
Base step. K[X)=0, then there are 2 cases.
Case 1: X=P. (T2oT1)[P]=T2[T1[P]]=T2[P]=P. Case 2: X=λ. (T2oT1)[λ]=T2[T1[λ]]=T2[λ]=λ.
Inductive step. K[X]≥1. As an inductive hypothesis we have that (T2oT1)[X1]=X1, (T2oT1)[X2]=X2.
There are 4 cases. Case 3. X=−X1. (T2oT1) [X1]=T2[T1[X1]]=T2[−{T1[X1]−}]=−(T2oT1)[X1]=−X1.
Case 4. X=X1∩X2. (T2oT1)[X1∩X2]=T2[T1[X1∩X2]]=T2[T1[X1]T1[X2]]=T2[T1[X1]] ∩T2[T1[X2]]
=X1∩X2.
Case 5. X=X1⊃X2. (T2oT1)[X1⊃X2]=T2[T1[X1⊃X2]]=T2[(T1[X1](T1[X2]))]
= T2[T1[X1]]⊃T2[T1[X2]]=X1⊃ X2.
Case 6. X=X1∪X2. (T2oT1)[X1∪X2]=T2[T1[X1∪X2]]=T2[((T1[X1])(T1[X2]))]
= T2[T1[X1]]∪T2[T1[X2]] = X1∪ X2.
By the principle of mathematical induction it has been proved that (T2oT1)=$Id_{FG.}$
Parts c and d. Direct consequence of parts a and b.

**Proposition 38.** For G,H∈GET, and X,Y∈FT. a) G∈TGET if and only if G"∈TT. b) X'∈TT yes and only if X∈TGET. c) G>>H if and only if H" is *consequence of* G". d) X'>>Y' if and only if Y it's *consequence of* X.

Proof part a. By proposition 36 we have that, if G∈TGET then G"∈TT, furthermore, by proposition 32 we have that, if G"∈TGET then (G")'∈TGET, but by proposition 37 we know that, (G")'=G, resulting that, if G"∈TGET then G∈TGET, and since we have the reciprocal, we conclude that, G∈TGET if and only if G"∈TT.
Proof part b. By proposition 30 we have that, if X∈TT then X'∈TGET, furthermore, by proposition 37 we have that, if X'∈TGET then (X')"∈TT, but by proposition 37 we know that, (X')"=X, resulting that, if X'∈TGET then X∈TT, and since we have the reciprocal, we conclude that, X'∈TGET if and only if X∈TT.
Proof part c. By proposition 30 we have, if X" implies Y" then [X"]'>>[Y"]', by proposition 37 we have [X"]'=X and [Y"]'=Y, so if X" implies Y" then X>>Y, in addition by proposition 36 we have the reciprocal. Therefore, G>>H if and only if G" implies H".
Proof part d. by proposition 38 we have, if X'⟩⟩Y' then [X']">>[Y"]', by proposition 37 we have [X']"=X and [Y']"=Y, so if X' implies Y' then X'>>Y', by proposition 30 we have the reciprocal. Therefore, X'>>Y' if and only if X' implies Y'.

## VII. Conclusions

**Definition 20.** For P∈FA; X,Y∈FT. For SD a deductive system with a negation operator N and for X a formula for SD. SD is said to be *paraconsistent* when SD does not derive all SD formulas from X and NX.

**Conclusion 1.** For G,H,K∈GET. a) G"⊃(−G"⊃H")∉VT. b) GET is paraconsistent. c) GT is paraconsistent.

Proof part a. Consider the models with 2 possible worlds M and N, such that M is the actual world, M<N, M<M, and N<N, where M(G")=1, M(H")=0, and N(G")=0. Since N(G")=0 and M<N, by V we infer M(−G")=1, and we have M(H")=0, applying V⊃ we derive M(−G''⊃H'')=0, and since M(G")=1, using V again⊃ we deduce M(G''⊃(−G''⊃H''))=0. Therefore, G"⊃(−G"⊃H")∉VT.
Proof part b. Applying proposition 38 yields (G(({G}(H))))∉TGET, which by Rules R3 and R6 implies that, this is not the case: G{G}>>H. Therefore, GET is paraconsistent.
Proof part c. Applying proposition 38 in part b.

**Conclusion 2.** For Mo=(S, <, Ma, V) a GT model.
Mo(+X'⊃[−(+X'∩Y')]⊃−Y')=1 equivalent to RT: (∀N,M,P∈S)(M<N y N<P implies M<P).

Proof. Suppose that Mo(+X'⊃[−(+X'∩Y')]⊃−Y')=0 means that Ma(+X'⊃[−(+X'∩Y')]⊃−Y')=0, by V⊃ result Ma(+X')=1, Ma(−(+X'∩Y'))=1 y Ma(−Y')=0, using V− follows the existence of P∈S, Ma<P y P(+X'∩Y')=0, also for all N∈S, Ma<N implies N(Y')=1, and since Ma<P then P(Y')=1, by V∩ is derived P(+X')=0, which by V− implies the existence of Q∈S, P<Q, Q(X')=0, as Ma(+X')=1, according to V⊃ it follows that for all N∈S, Ma<N implies N(X')=1, as Ma<P then P(X')=1.
If it is not satisfied (∀N,M,P∈S)(M<N y N<P implies M<P) then Mo(+X'⊃[−(+X'∩Y')]⊃−Y')=0, since Ma<P, Ma<Ma, P<P, P<Q, Q<Q, Q(X')=0, P(X')=1 and P(Y')=1, but it is not the case that M<Q.
If it is satisfied (∀N,M,P∈S)(M<N y N<P implies M<P) then, since Ma<P and P<Q it follows that Ma<Q, and since Ma(+X')=1, according to V⊃ it follows that Q(X')=1, which is not the case, hence Mo(+X'⊃[−(+X'∩Y')]⊃−Y')=1.

**Conclusion 3.** For X,Y∈FT. a) +X⊃[−Y⊃−(+X∩Y)]∈TT. b) +X∩−Y≡+X∩−(+X∩Y)]∈TT equivalent to *Y'{X'}⇔*Y'{*Y'X'} in TGET. c) *Y'{X'}⇔*Y'{*Y'X'} in GET implies strong I-D: *Y'[X']$_{1cq}$⇔*Y'[*Y'X']$_{1cq}$ in GET.

Proof part a. By the application of rule R2, (+X⊃[−Y⊃−(+X∩Y)])'∈TGET, following proposition 38, it is concluded that +X⊃[−Y⊃−(+X∩Y)]∈TT.
Proof part b. It is derived from proposition 38.
Proof part c. Direct consequence of proposition 34.

Observation. The test in conclusion 2 shows that (+X'⊃[−(+X'∩Y')]⊃−Y'∉TT, so a new deductive system and a new graphing system are required to properly utilize conclusion 3.

**Definition 26.** The GET4 system of existential graphs is defined as the result of adding the strong I-D rule to GET. The GT4 deductive system is defined as the result of adding the formula +X⊃[−(+X∩Y)⊃−Y] as an axiom to GT.

**Conclusion 4.** For G,H∈GET4 and X,Y∈FT4. a) G∈TGET4 if and only if G"∈TT4. b) X'∈TT4 if and only if X∈TGET4. c) G>>H if and only if H" is *Consequence of* G". d) X'>>Y' if and only if Y it's *Consequence of* X.

Proof. Direct consequence of propositions 30, 35 to 38, and conclusions 2 and 3.

**Conclusion 5.** For G,H,K∈GET4. a) G"⊃(−G"⊃H")∉VT4. b) GET4 is paraconsistent. c) GT4 is paraconsistent.

Proof. The same reasoning as the conclusion 1 proof.

**Conclusion 6.** Clearly, GET4 matches Zeman's Gamma-4. In addition, in Gamma-4.2 Gamma-5 is valid conclusion 5, consequently, they are also paraconsistent.